\documentclass[11pt]{amsart}

\usepackage{tikz,amsmath,amssymb}

\theoremstyle{theorem}
\newtheorem{theorem}{Theorem}[section]

\newtheorem{corollary}[theorem]{Corollary}

\theoremstyle{definition}
\newtheorem{remark}[theorem]{Remark}
\newtheorem{definition}[theorem]{Definition}

\newtheorem{example}[theorem]{Example}
\newtheorem{problem}[theorem]{Problem}

\newcommand{\R}{\mathbb{R}}
\newcommand{\Z}{\mathbb{Z}}

\newcommand{\bv}{\mathbf{v}}
\newcommand{\LP}{\mathbf{L}}
\newcommand{\bw}{\mathbf{w}}
\newcommand{\bov}{\mathbf{v}}
\newcommand{\bx}{\mathbf{x}}
\newcommand{\ba}{\mathbf{a}}
\newcommand{\bo}{\mathbf{0}}
\newcommand{\vol}{\mathrm{vol}}
\newcommand{\Vol}{\mathrm{Vol}}

\newcommand\commentout[1]{}

\begin{document}



\title{A Brief Survey on Lattice Zonotopes}

\author{Benjamin Braun}
\address{Department of Mathematics\\
         University of Kentucky\\
         Lexington, KY 40506--0027\\
https://sites.google.com/view/braunmath/}
\email{benjamin.braun@uky.edu}

\author{Andr\'es R. Vindas-Mel\'endez}
\address{Department of Mathematics\\
         University of Kentucky\\
         Lexington, KY 40506--0027\\
https://math.as.uky.edu/users/arvi222}
\email{andres.vindas@uky.edu}

\date{14 August 2018}



\begin{abstract}
Zonotopes are a rich and fascinating family of polytopes, with connections to many areas of mathematics.
In this article we provide a brief survey of classical and recent results related to lattice zonotopes.
Our emphasis is on connections to combinatorics, both in the sense of enumeration (e.g. Ehrhart theory) and combinatorial structures (e.g. graphs and permutations).
\end{abstract}

\maketitle

\section{Introduction}\label{sec:intro}

Zonotopes are a rich and fascinating family of polytopes, with connections to many areas of mathematics.
In this article we provide a brief survey of classical and recent results related to \emph{lattice} zonotopes, i.e. Minkowski sums of line segments with endpoints in the lattice $\Z^n$.
Our emphasis is on connections to combinatorics, both in the sense of enumeration (e.g. Ehrhart theory) and combinatorial structures (e.g. graphs and permutations).
Our primary goal in this article is to maintain a level of presentation accessible to beginning graduate students working in algebraic and geometric combinatorics or related fields.
We make no effort here to be comprehensive; for example, we omit from our discussion the deep connections between zonotopes, hyperplane arrangements, and oriented matroids.
Our selection of topics is based purely on personal taste, and includes various topics we find interesting for one reason or another.
Our hope is that this survey will serve as a stepping stone to a deeper investigation of lattice zonotopes for interested readers.

\section{Zonotope Basics}\label{sec:basics}

Zonotopes can be defined using either Minkowski sums or projections of cubes.

\begin{definition}\label{def:msum}
Consider the polytopes, $P_1, P_2,  \dots, P_m \subset \R^n$. 
We define the \emph{Minkowski sum} of the $m$ polytopes as 
\begin{equation*}
P_1+P_2+\cdots+P_m:=\{x_1+x_2+\cdots+x_m:x_j \in P_j \text{ for } 1\leq j\leq m\} \, .
\end{equation*}
\end{definition}

Given $\bv,\bw\in \R^n$, we write $[\bv,\bw]$ for the line segment from $\bv$ to $\bw$.

\begin{example}
Consider the Minkowski sum of $[(0,0),(1,0)]$, $[(0,0),(0,1)]$, and $[(0,0),(1,1)]$.
The Minkowski sum of the first two segments is a unit square.
Taking the Minkowski sum of this square with the line segment $[(0,0),(1,1)]$ can be visualized as sliding the square up and to the right along the line segment, with the resulting polytope consisting of all points touched by the square during the sliding movement.
See Figure~\ref{fig:minksum}.

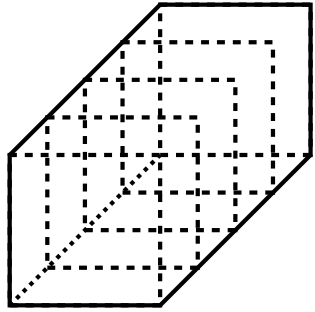
\begin{figure}
\begin{center}
\begin{tikzpicture}[scale=2]
\draw[black, ultra thick] (0,0) -- (1,0) -- (2,1) -- (2,2) -- (1,2) -- (0,1) -- cycle;
\draw[black, ultra thick, dashed] (.5,.5) -- (1.5,0.5) -- (1.5,1.5) -- (0.5,1.5) -- cycle;
\draw[black, ultra thick, dashed] (1,1) -- (2,1) -- (2,2) -- (1,2) -- cycle;
\draw[black, ultra thick, dashed] (.25,.25) -- (1.25,.25) -- (1.25,1.25) -- (.25,1.25) -- cycle;
\draw[black, ultra thick, dashed] (.75,.75) -- (1.75,.75) -- (1.75,1.75) -- (.75,1.75) -- cycle;
\draw[black, ultra thick, dashed] (0,0) -- (1,0) -- (1,1) -- (0,1) -- cycle;
\draw[black, ultra thick, dotted] (0,0) -- (1,1);
\end{tikzpicture}
\end{center}
\caption{The Minkowski sum of $[(0,0),(1,0)]$, $[(0,0),(0,1)]$, and $[(0,0),(1,1)]$.}
\label{fig:minksum}
\end{figure}
\end{example}

\begin{definition}\label{def:zonotope}
Consider $m$ vectors $\bv_1,\ldots,\bv_m$ in $\R^n$ and their corresponding line segments $[\bo,\bv_j]$.
The \emph{zonotope} corresponding to $\bv_1,\ldots,\bv_m$ is defined to be the Minkowski sum of these line segments: 
\begin{equation*}
Z(\bv_1,\bv_2,\dots, \bv_m):=\{\lambda_1\bv_1+\lambda_2\bv_2+\cdots+\lambda_m\bv_m: 0\leq \lambda_j\leq 1\}= [\bo,\bv_1]+\cdots+[\bo,\bv_m]\, .
\end{equation*}
We call any polytope that is translation-equivalent to an object of this type a \emph{zonotope}.
When each $\bv_j\in \Z^n$, we say this is a \emph{lattice zonotope}.
\end{definition}

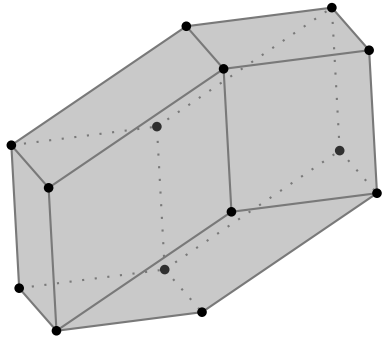
\begin{figure}
\begin{center}
\begin{tikzpicture}%
	[x={(0.247809cm, -0.283385cm)},
	y={(0.967407cm, 0.123809cm)},
	z={(-0.052102cm, 0.950981cm)},
	scale=2.000000,
	back/.style={loosely dotted, thick},
	edge/.style={color=gray!95!black, thick},
	facet/.style={fill=gray!95!black,fill opacity=0.400000},
	vertex/.style={inner sep=1pt,circle,draw=black!25!black,fill=black!75!black,thick,anchor=base}]
%
%
\coordinate (0.00000, 0.00000, 0.00000) at (0.00000, 0.00000, 0.00000);
\coordinate (0.00000, 0.00000, 1.00000) at (0.00000, 0.00000, 1.00000);
\coordinate (0.00000, 1.00000, 0.00000) at (0.00000, 1.00000, 0.00000);
\coordinate (0.00000, 1.00000, 1.00000) at (0.00000, 1.00000, 1.00000);
\coordinate (1.00000, 0.00000, 0.00000) at (1.00000, 0.00000, 0.00000);
\coordinate (1.00000, 0.00000, 1.00000) at (1.00000, 0.00000, 1.00000);
\coordinate (1.00000, 1.00000, 0.00000) at (1.00000, 1.00000, 0.00000);
\coordinate (2.00000, 2.00000, 2.00000) at (2.00000, 2.00000, 2.00000);
\coordinate (1.00000, 1.00000, 2.00000) at (1.00000, 1.00000, 2.00000);
\coordinate (1.00000, 2.00000, 1.00000) at (1.00000, 2.00000, 1.00000);
\coordinate (1.00000, 2.00000, 2.00000) at (1.00000, 2.00000, 2.00000);
\coordinate (2.00000, 1.00000, 1.00000) at (2.00000, 1.00000, 1.00000);
\coordinate (2.00000, 1.00000, 2.00000) at (2.00000, 1.00000, 2.00000);
\coordinate (2.00000, 2.00000, 1.00000) at (2.00000, 2.00000, 1.00000);
\draw[edge,back] (0.00000, 0.00000, 0.00000) -- (0.00000, 1.00000, 0.00000);
\draw[edge,back] (0.00000, 0.00000, 1.00000) -- (0.00000, 1.00000, 1.00000);
\draw[edge,back] (0.00000, 1.00000, 0.00000) -- (0.00000, 1.00000, 1.00000);
\draw[edge,back] (0.00000, 1.00000, 0.00000) -- (1.00000, 1.00000, 0.00000);
\draw[edge,back] (0.00000, 1.00000, 0.00000) -- (1.00000, 2.00000, 1.00000);
\draw[edge,back] (0.00000, 1.00000, 1.00000) -- (1.00000, 2.00000, 2.00000);
\draw[edge,back] (1.00000, 2.00000, 1.00000) -- (1.00000, 2.00000, 2.00000);
\draw[edge,back] (1.00000, 2.00000, 1.00000) -- (2.00000, 2.00000, 1.00000);
\node[vertex] at (1.00000, 2.00000, 1.00000)     {};
\node[vertex] at (0.00000, 1.00000, 0.00000)     {};
\node[vertex] at (0.00000, 1.00000, 1.00000)     {};
\fill[facet] (1.00000, 1.00000, 2.00000) -- (2.00000, 1.00000, 2.00000) -- (2.00000, 2.00000, 2.00000) -- (1.00000, 2.00000, 2.00000) -- cycle {};
\fill[facet] (2.00000, 1.00000, 1.00000) -- (2.00000, 2.00000, 1.00000) -- (2.00000, 2.00000, 2.00000) -- (2.00000, 1.00000, 2.00000) -- cycle {};
\fill[facet] (1.00000, 0.00000, 1.00000) -- (0.00000, 0.00000, 1.00000) -- (0.00000, 0.00000, 0.00000) -- (1.00000, 0.00000, 0.00000) -- cycle {};
\fill[facet] (2.00000, 1.00000, 2.00000) -- (1.00000, 0.00000, 1.00000) -- (0.00000, 0.00000, 1.00000) -- (1.00000, 1.00000, 2.00000) -- cycle {};
\fill[facet] (2.00000, 1.00000, 2.00000) -- (1.00000, 0.00000, 1.00000) -- (1.00000, 0.00000, 0.00000) -- (2.00000, 1.00000, 1.00000) -- cycle {};
\fill[facet] (2.00000, 2.00000, 1.00000) -- (1.00000, 1.00000, 0.00000) -- (1.00000, 0.00000, 0.00000) -- (2.00000, 1.00000, 1.00000) -- cycle {};
\draw[edge] (0.00000, 0.00000, 0.00000) -- (0.00000, 0.00000, 1.00000);
\draw[edge] (0.00000, 0.00000, 0.00000) -- (1.00000, 0.00000, 0.00000);
\draw[edge] (0.00000, 0.00000, 1.00000) -- (1.00000, 0.00000, 1.00000);
\draw[edge] (0.00000, 0.00000, 1.00000) -- (1.00000, 1.00000, 2.00000);
\draw[edge] (1.00000, 0.00000, 0.00000) -- (1.00000, 0.00000, 1.00000);
\draw[edge] (1.00000, 0.00000, 0.00000) -- (1.00000, 1.00000, 0.00000);
\draw[edge] (1.00000, 0.00000, 0.00000) -- (2.00000, 1.00000, 1.00000);
\draw[edge] (1.00000, 0.00000, 1.00000) -- (2.00000, 1.00000, 2.00000);
\draw[edge] (1.00000, 1.00000, 0.00000) -- (2.00000, 2.00000, 1.00000);
\draw[edge] (2.00000, 2.00000, 2.00000) -- (1.00000, 2.00000, 2.00000);
\draw[edge] (2.00000, 2.00000, 2.00000) -- (2.00000, 1.00000, 2.00000);
\draw[edge] (2.00000, 2.00000, 2.00000) -- (2.00000, 2.00000, 1.00000);
\draw[edge] (1.00000, 1.00000, 2.00000) -- (1.00000, 2.00000, 2.00000);
\draw[edge] (1.00000, 1.00000, 2.00000) -- (2.00000, 1.00000, 2.00000);
\draw[edge] (2.00000, 1.00000, 1.00000) -- (2.00000, 1.00000, 2.00000);
\draw[edge] (2.00000, 1.00000, 1.00000) -- (2.00000, 2.00000, 1.00000);
\node[vertex] at (0.00000, 0.00000, 0.00000)     {};
\node[vertex] at (0.00000, 0.00000, 1.00000)     {};
\node[vertex] at (1.00000, 0.00000, 0.00000)     {};
\node[vertex] at (1.00000, 0.00000, 1.00000)     {};
\node[vertex] at (1.00000, 1.00000, 0.00000)     {};
\node[vertex] at (2.00000, 2.00000, 2.00000)     {};
\node[vertex] at (1.00000, 1.00000, 2.00000)     {};
\node[vertex] at (1.00000, 2.00000, 2.00000)     {};
\node[vertex] at (2.00000, 1.00000, 1.00000)     {};
\node[vertex] at (2.00000, 1.00000, 2.00000)     {};
\node[vertex] at (2.00000, 2.00000, 1.00000)     {};
\end{tikzpicture}
\end{center}
\caption{$Z((1,0,0),(0,1,0),(0,0,1),(1,1,1))$}
\label{fig:3dminksum} 
\end{figure}

\begin{example}
The Minkowski sum in Figure~\ref{fig:minksum} is $Z((1,0),(0,1),(1,1))$.
The Minkowski sum in Figure~\ref{fig:3dminksum} is $Z((1,0,0),(0,1,0),(0,0,1),(1,1,1))$.
\end{example}

\begin{definition}\label{def:centsymzonotope}
The zonotope $Z_0(\bv_1,\bv_2,\dots,\bv_m):=Z(\pm\bv_1,\pm\bv_2,\dots,\pm\bv_m)$ is symmetric about the origin, that is, it has the property that $\bx \in Z_0$ if and only if $-\bx \in Z_0$; we call $Z_0$ a \emph{centrally symmetric zonotope} defined by $\bv_1,\ldots,\bv_m$.
Note that $Z_0$ can be obtained as $Z(2\bv_1,2\bv_2,\dots,2\bv_m)-(\bv_1+\cdots+\bv_m)$.
\end{definition}

An alternative definition of a zonotope is as a projection of the unit cube.
In this case, let $A$ denote the matrix with columns given by $\bv_1,\ldots,\bv_m$.
Then by definition it is immediate that $Z(\bv_1,\ldots,\bv_m)$ is equal to $A\cdot [0,1]^m$.
In some circumstances it is more convenient to work with this projection-based definition of zonotopes.

\begin{remark}
Zonotopes are deceptively simple to define, yet even the most elementary zonotopes are mathematically rich.
For example, the unit cube $[0,1]^n$ is itself a zonotope, and the survey paper by Zong~\cite{Zong} shows that the mathematical properties of this object are both broad and deep.
\end{remark}

In polyhedral geometry, it is typically useful when objects can be decomposed as unions of simpler objects, e.g. the theory of subdivisions and triangulations.
Zonotopes admit a particularly nice decomposition into parallelepipeds; parts of the boundaries of these parallelepipeds can be removed resulting in a disjoint decomposition.
More precisely, suppose that $\bw_1,\bw_2,\dots,\bw_k\in \R^n$ are linearly independent, and let $\sigma_1,\sigma_2,\dots,\sigma_k\in \{\pm1\}$. 
Then we define 
\begin{equation*}
\Pi^{\sigma_1,\sigma_2,\dots,\sigma_k}_{\bw_1,\bw_2,\dots,\bw_k}:=\left\{\lambda_1\bw_1+\lambda_2\bw_2+\cdots+\lambda_k\bw_k: \begin{aligned} 0\leq \lambda_j< 1 \text{ if } \sigma_j&=-1 \\
 0< \lambda_j\leq 1 \text{ if } \sigma_j&=1\end{aligned} \right\}
\end{equation*}
to be the \emph{half-open parallelepiped} generated by $\bw_1,\bw_2,\dots, \bw_k$. 
The signs $\sigma_1,\sigma_2,\dots, \sigma_k$ keep track of the facets of the parallelepiped that are either included or excluded from the closure of the parallelepiped. 

\begin{theorem}[Shephard~\cite{Shephard}, Theorem 54]
\label{thm:decomp}
The zonotope $Z(\bv_1,\bv_2,\dots, \bv_m)$ can be written as a disjoint union of translates of $\Pi^{\sigma_1,\sigma_2,\dots,\sigma_k}_{\bw_1,\bw_2,\dots,\bw_k}$, where $\{\bw_1,\bw_2,\dots,\bw_k\}$ ranges over all linearly independent subsets of $\{\bv_1,\bv_2,\dots,\bv_m\}$, each equipped with an appropriate choice of signs $\sigma_1,\sigma_2,\dots,\sigma_k$. 
\end{theorem}

\begin{figure}
\begin{center}
\begin{tikzpicture}%
[	back/.style={loosely dotted, thick},
	edge/.style={color=gray!95!black, thick},
	facet/.style={fill=gray!95!black,fill opacity=0.400000},
	mydot/.style={circle,inner sep=1pt,circle,draw=black!25!black,fill=white!75!white,thick,anchor=base},
  	vertex/.style={inner sep=1pt,circle,draw=black!25!black,fill=black!75!black,thick,anchor=base}]
%
%
%
\coordinate (A) at (0.00000, 0.00000);
\coordinate (B) at (0.0000, 0.30000);
\coordinate (C) at (0.00000, 2.30000);

\node[vertex] at (0,0)     {};
\node[mydot] at (0.0000, 0.30000) {};
\draw (B) to (C);

\coordinate (D) at (2,0);
\coordinate (E) at (2, .3); 
\coordinate (F) at (2,2.3);

\draw (F) to (E); 
\node[vertex] at (2,0)     {};
\node[mydot] at (2.0000, 0.30000) {};

\coordinate (G) at (2.3,.3);
\coordinate (H) at (4.3,2.3);
\draw (G) to (H);
\node[mydot] at (2.3000, 0.30000) {};

\coordinate (I) at (2.3,.6); 
\coordinate (J) at (2.3, 2.6);
\draw[dashed] (I) to (J);
\node[mydot] at (2.3, 2.6) {};

\coordinate (K) at (4.3,2.6); 
\draw[dashed] (I) to (K);
\node[mydot] at (4.3,2.6) {};

\coordinate (L) at (4.3,4.6);
\draw (J) to (L);
\draw (K) to (L);
\draw[gray,top color=gray, bottom color=gray, fill opacity=0.25] (J) -- (L) -- (K) -- (I) -- cycle;

\coordinate (M) at (6,0);
\coordinate (N) at (6, .3); 
\coordinate (O) at (6,2.3);

\draw (O) to (N); 
\node[vertex] at (6,0)     {};
\node[mydot] at (6.0000, 0.30000) {};

\coordinate (P) at (6.3,.3);
\coordinate (Q) at (8.3,2.3);
\draw (P) to (Q);
\node[mydot] at (6.3000, 0.30000) {};

\coordinate (R) at (6.3,.6); 
\coordinate (S) at (6.3, 2.6);
\draw[dashed] (R) to (S);
\node[mydot] at (6.3, 2.6) {};

\coordinate (T) at (8.3,2.6); 
\draw[dashed] (R) to (T);
\node[mydot] at (8.3,2.6) {};

\coordinate (U) at (8.3,4.6);
\draw (S) to (U);
\draw (T) to (U);
\draw[gray,top color=gray, bottom color=gray, fill opacity=0.25] (S) -- (U) -- (T) -- (R) -- cycle;

\coordinate (V) at (6.3,.15); 
\coordinate (W) at (8.3,.15);
\draw (V) to (W); 
\node[mydot] at (6.3,.15) {};

\coordinate (X) at (6.5, .25);
\coordinate (Y) at (8.5, .25); 
\coordinate (Z) at (8.6,2.3);
\coordinate (a) at (10.6, 2.3);
\draw[dashed] (X) to (Y);
\draw[dashed](X) to (Z);
\draw(Z) to (a);
\draw(Y) to (a);
\node[mydot] at (Z) {};
\node[mydot] at (Y) {};
\draw[gray,top color=gray, bottom color=gray, fill opacity=0.25] (X) -- (Y) -- (a) -- (Z) -- cycle;

\coordinate (b) at (8.6, 2.6);
\coordinate (c) at (10.6, 2.6);
\coordinate (d) at (8.6, 4.6); 
\coordinate (e) at (10.6,4.6);
\draw[dashed] (b) to (c);
\draw[dashed] (b) to (d);
\draw (e) to (d); 
\draw[gray,top color=gray, bottom color=gray, fill opacity=0.25] (b) -- (c) -- (e) -- (d) -- cycle;

\end{tikzpicture}
\end{center}
\caption{A zonotopal decomposition of $Z((0,1),(1,0),(1,1)$}
\label{fig:zonotopedecomp} 
\end{figure}
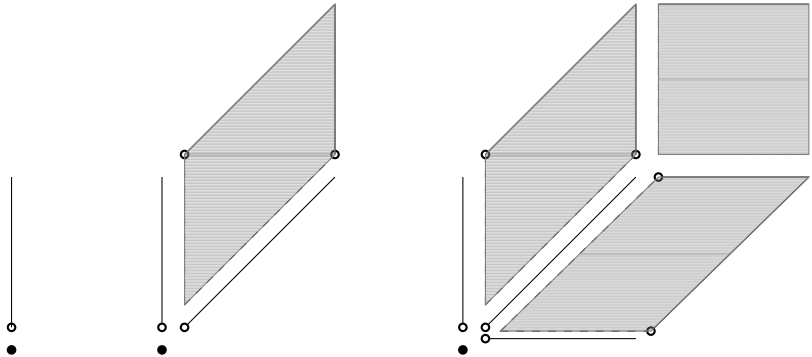

For a complete proof of Theorem~\ref{thm:decomp}, see Lemma 9.1 of the textbook by Beck and Robins~\cite{BeckRobinsCCD}.  
Figure \ref{fig:zonotopedecomp} is an illustration of the decomposition of the zonotope $Z((0,1),(1,0),(1,1))\subset \Z^2$ as suggested by Shephard's theorem. 


\section{Ehrhart Theory}\label{sec:ehrhartbackground}

A \emph{lattice polytope} $P\subset \R^n$ of dimension $d$ is the convex hull of finitely many points in $\Z^n$ that together affinely span a $d$-dimensional hyperplane.
For $t\in\Z_{>0}$, set $tP:=\{tp:p\in P\}$, and let $L_P(t)=|\Z^n\cap tP|$.  
Ehrhart~\cite{Ehrhart} proved a statement equivalent to the following.
Recall that the set 
\[
\left\{{t+d-i\choose d}:i=0,1,\ldots,d\right\}
\]
is a basis for polynomials of degree $d$, where 
\[
{t+d-i\choose d}=(1/d!)(t+d-i)(t+d-i-1)(t+d-i-2)\cdots(t-i+1)
\]
is clearly a polynomial in $t$ of degree $d$.
For any lattice polytope $P$, there exist rational values $c_0,c_1,\ldots,c_d$ and $h_0^*,h_1^*,\ldots,h_d^*$ such that
\begin{equation*}
L_P(t)=\sum_{i=0}^dh_i^*{t+d-i\choose d} = \sum_{i=0}^dc_it^i \, .
\end{equation*}
The polynomial $L_P(t)$ is called the \emph{Ehrhart polynomial of $P$} and has connections to commutative algebra, algebraic geometry, combinatorics, and discrete and convex geometry.  
Stanley~\cite{StanleyDecompositions} proved that $h_i^*\in \Z_{\geq 0}$ for all $i$.
We call the polynomial $h^\ast(P;x):=h_0^\ast+h_1^\ast x+\cdots+h_d^\ast x^d$  encoding the $h^*$-coefficients the \emph{$h^\ast$-polynomial} (or \emph{$\delta$-polynomial}) of $P$.
The coefficients of $h^*(P;x)$ form the \emph{$h^*$-vector of $P$}.

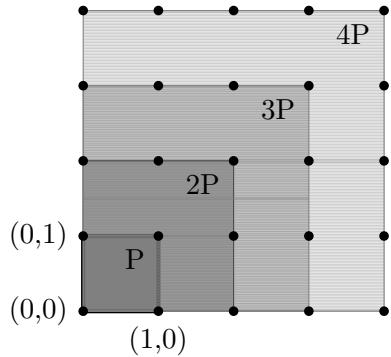
\begin{figure}
\begin{center}
\begin{tikzpicture}
[back/.style={loosely dotted, thick},
	edge/.style={color=gray!95!black, thick},
	facet/.style={fill=gray!95!black,fill opacity=0.400000},
	mydot/.style={circle,inner sep=1pt,circle,draw=black!25!black,fill=white!75!white,thick,anchor=base},
  	vertex/.style={inner sep=1pt,circle,draw=black!25!black,fill=black!75!black,thick,anchor=base}]
	
\draw[black, ultra thick][black,top color=gray, bottom color=gray, fill opacity=1] (0,0) -- (1,0) -- (1,1) -- (0,1) -- cycle;
\draw[black][black,top color=gray, bottom color=gray, fill opacity=0.5] (0,0) -- (2,0) -- (2,2) -- (0,2) -- cycle;
\draw[black][gray,top color=gray, bottom color=gray, fill opacity=0.3]  (0,0) -- (3,0) -- (3,3) -- (0,3) -- cycle; 
\draw[black][gray,top color=gray, bottom color=gray, fill opacity=0.2]  (0,0) -- (4,0) -- (4,4) -- (0,4) -- cycle;
\node[vertex, label=left:{(0,0)}] at (0,0) {};
\node[vertex, label=below:{(1,0)}]  at (1,0) {};
\node[vertex, label=left:{(0,1)}] at (0,1) {};
\node[vertex, label=below left:{P}] at (1,1) {};
\node[vertex] at (2,0) {};
\node[vertex] at (2,1) {};
\node[vertex, label=below left:{2P}] at (2,2) {};
\node[vertex] at (2,3) {};
\node[vertex] at (2,4) {};
\node[vertex] at (0,2) {};
\node[vertex] at (1,2) {};
\node[vertex] at (1,3) {};
\node[vertex] at (1,4) {};
\node[vertex] at (0,3) {};
\node[vertex] at (0,4) {};
\node[vertex] at (3,0) {};
\node[vertex] at (3,1) {};
\node[vertex] at (3,2) {};
\node[vertex, label=below left:{3P}] at (3,3) {};
\node[vertex] at (3,4) {};
\node[vertex] at (4,0) {};
\node[vertex] at (4,1) {};
\node[vertex] at (4,2) {};
\node[vertex] at (4,3) {};
\node[vertex, label=below left:{4P}] at (4,4) {};

\end{tikzpicture}
\end{center}
\caption{The unit square $P=[0,1]^2\subset \Z^2$ and some of its dilates.}
\label{fig:unitsquare}
\end{figure}

\begin{example}
The unit square as depicted in Figure \ref{fig:unitsquare} has Ehrhart polynomial $L_{[0,1]^2}(t)=t^2+2t+1$ and $h^*([0,1]^2;x)=1+x$.
\end{example}

Various properties of $P$ are reflected in its $h^*$-polynomial.
For example, $\vol(P)=(\sum_ih_i^*)/d!$, where $\vol(P)$ denotes the Euclidean volume of $P$ with respect to the integer lattice contained in the hyperplane spanned by $P$.  
We therefore define the \emph{normalized volume} of $P$ to be $\Vol(P)=\sum_ih^*_i$.
Further, it is known that $h_0^*=1$ for all $P$, and that $h_d^*$ is equal to the number of lattice points in the relative (topological) interior of $P$ within the affine span of $P$.
Another interesting combinatorial property displayed by  $h^\ast(P;x)$ for some lattice polytopes is \emph{unimodality}.  
A polynomial $a_0+a_1x+\cdots+a_dx^d$ is called \emph{unimodal} if there exists an index $j$, $0\leq j\leq d$, such that $a_{i-1}\leq a_i$ for $i\leq j$, and $a_i\geq a_{i+1}$ for $i\geq j$.  
Unimodality of $h^*$-polynomials is an area of active research ~\cite{BraunUnimodalSurvey}.


\section{Ehrhart Polynomials for Zonotopes}\label{sec:ehrhartpoly}

For a lattice zonotope $Z$, Stanley proved the following description of the coefficients of $L_Z(t)$.

\begin{theorem}[Stanley~\cite{stanleyzonotope}, Theorem 2.2]\label{thm:ehrpoly}
Let $Z:=Z(\bv_1,\dots,\bv_m)$ be a zonotope generated by the integer vectors $\bv_1,\dots,\bv_m$. 
Then the Ehrhart polynomial of $Z$ is given by 
\begin{equation*}
L_Z(t)=\sum_{S}m(S)t^{|S|} \, ,
\end{equation*}
where $S$ ranges over all linearly independent subsets of $\{\bv_1,\dots,\bv_m\}$, and $m(S)$ is the greatest common divisor of all minors of size $|S|$ of the matrix whose columns are the elements of $S$. 
\end{theorem}

The proof of Theorem~\ref{thm:ehrpoly} relies on Theorem~\ref{thm:decomp}, and is our first example of the usefulness of half-open decompositions of zonotopes.
A generalization of Stanley's theorem was recently given by Hopkins and Postnikov~\cite{HopkinsPostnikovEhrhart}.

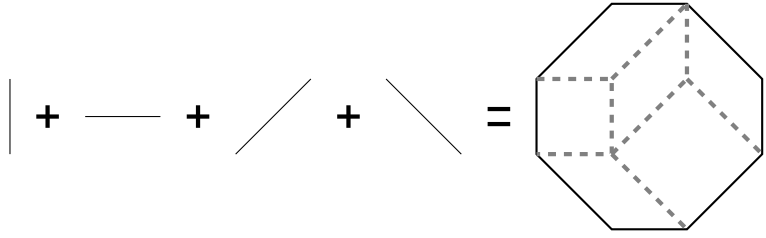
\begin{figure}
\begin{center}
\begin{tikzpicture}[scale=1]
\draw[black] (0,0) -- (0,1);
\draw[black, ultra thick] (.5,.35) -- (.5,.65);
\draw[black, ultra thick] (.35, .50) -- (.65,.50);
\draw[black] (1,.5) -- (2,.5);
\draw[black, ultra thick] (2.5,.35) -- (2.5,.65);
\draw[black, ultra thick] (2.35, .50) -- (2.65,.50);
\draw[black] (3,0) -- (4,1);
\draw[black, ultra thick] (4.5,.35) -- (4.5,.65);
\draw[black, ultra thick] (4.35, .50) -- (4.65,.50);
\draw[black] (5,1) -- (6,0);
\draw[black, ultra thick] (6.35,.40) -- (6.65,.40);
\draw[black, ultra thick] (6.35,.60) -- (6.65,.60);
\draw[black, thick] (7,0) -- (7,1) -- (8,2) -- (9,2) -- (10,1) -- (10,0) -- (9,-1) -- (8, -1) -- cycle; 
\draw[gray, ultra thick, dashed] (7,1) -- (8,1) -- (8,0) -- (7,0);
\draw[gray, ultra thick, dashed] (8,1) -- (9,2);
\draw[gray, ultra thick, dashed] (8,0) -- (9,1)--(9,2);
\draw[gray, ultra thick, dashed] (9,1) -- (10,0);
\draw[gray, ultra thick, dashed] (8,0) -- (9,-1);
\end{tikzpicture}
\end{center}
\caption{The zonotope $Z((0,1),(1,0),(1,1),(1,-1))$ is a Minkowski sum of four line segments.}
\label{fig:fourlinesegmentzonotope}
\end{figure}

\begin{example}\label{ex:stanleytheorem}
Consider the zonotope $Z((0,1),(1,0),(1,1),(1,-1))\subset \Z^2$ depicted in Figure \ref{fig:fourlinesegmentzonotope}. 
We compute
\begin{align*}
L_{Z((0,1),(1,0),(1,1),(1,-1))}(t) =&
\left|\mathrm{det}\begin{pmatrix} 
0 & 1 \\
1 & 0 
\end{pmatrix}\right|t^2  +\left|\mathrm{det}\begin{pmatrix} 
0 & 1 \\
1 & 1 
\end{pmatrix}\right|t^2+\left|\mathrm{det}\begin{pmatrix} 
0 & 1 \\
1 & -1
\end{pmatrix}\right|t^2 +
\\
& \left|\mathrm{det}\begin{pmatrix} 
1 & 1 \\
0 & 1 
\end{pmatrix}\right|t^2+\left|\mathrm{det}\begin{pmatrix} 
1 & 1 \\
0 & -1 
\end{pmatrix}\right|t^2+\left|\mathrm{det}\begin{pmatrix} 
1 & 1 \\
1 & -1
\end{pmatrix}\right|t^2 + \\
& \, \mathrm{gcd}(0,1)t+\mathrm{gcd}(1,0)t+\mathrm{gcd}(1,1)t+\mathrm{gcd}(1,-1)t+1\\
= & \, 7t^2+4t+1
\end{align*}
\end{example}

Recent work by D'Adderio and Moci has established a connection between $L_Z(t)$ and the arithmetic Tutte polynomial defined as follows by Moci~\cite{mociarithmetic}.

\begin{definition}\label{def:arithmetictutte}	
Consider a collection $A\subseteq \R^n$ of integer vectors. 
The \emph{arithmetic Tutte polynomial} is 
\begin{equation*}
M_A(x,y):=\sum_{B\subseteq A}m(B)(x-1)^{r(A)-r(B)}(y-1)^{|B|-r(B)} \, ,
\end{equation*}
where for each $B\subseteq A$, the multiplicity $m(B)$ is the index of $\Z B$ as a sublattice of $\mathrm{span}(B)\cap \Z^n$ and $r(A)$ and $r(B)$ are the size of the largest independent subset of $A$ and $B$, respectively.
If we use the vectors in $B$ as the columns of a matrix, then $m(B)$ equals the greatest common divisor of the minors of full rank. 
\end{definition}

The result of D'Adderio and Moci is to obtain the Ehrhart polynomial of $Z$ as a specialization of the arithmetic Tutte polynomial associated to the lattice points generating $Z$.

\begin{theorem}[D'Adderio and Moci~\cite{dadderiomociehrhart}, Theorem 3.2]\label{thm:ehrharttutte}
Let $Z:=Z(\bv_1,\dots,\bv_m)$ and let the columns of $A$ be given by $\bv_1,\dots,\bv_m$. Then 
\begin{equation*}
L_Z(t)=t^{r(A)}M_A\left(1+\frac{1}{t},1\right)\, . 
\end{equation*}
\end{theorem}

\begin{example}
Consider the zonotope $Z((0,1),(1,0),(1,1),(1,-1))$ from Example \ref{ex:stanleytheorem}. 
Then 
\begin{align*}
M_Z(x,y)&=(x-1)^2+(1+1+1+1)(x-1)+(1+1+1+1+1+2)+(y-1)\\
&=x^2+2x+3+y.
\end{align*}
The corresponding Ehrhart polynomial is 
\[L_Z(t)=t^2M_Z\left(1+\frac{1}{t},1\right)=7t^2+4t+1,
\] 
which agrees with Example \ref{ex:stanleytheorem}.
\end{example} 

From basic properties of Ehrhart polynomials, e.g. Ehrhart-MacDonald reciprocity, we obtain the following corollary observed by D'Adderio and Moci~\cite{dadderiomociehrhart}.

\begin{corollary}\label{cor:tuttecorollary}
\begin{enumerate}
\item $L_{Z^\circ}(t)=(-t)^{r(A)}M_A\left(1-\frac{1}{t},1\right)$ .
\item The volume of the zonotope $Z(A)$ is $M_A(1,1)$.
\item The zonotope $Z(A)$ contains $M_A(2,1)$ lattice points.
\item The zonotope $Z(A)$ contains $M_A(0,1)$ interior lattice points. 
\end{enumerate}
\end{corollary}

These techniques have been applied in interesting ways; for example, Ardila, Castillo, and Henley~\cite{ardilacastillohenley} have computed the arithmetic Tutte polynomials and Ehrhart polynomials for zonotopes defined by the classical root systems.

\section{Ehrhart $h^*$-Polynomials for Zonotopes}\label{sec:ehrharthstar}

Determining the $h^*$-polynomial for a lattice zonotope $Z$ is less straightforward computationally than determining $L_Z(t)$, yet the $h^*$-vector of $Z$ is of great interest. 
Two properties of $h^*$-polynomials that have been studied in recent research in Ehrhart theory are unimodality and real-rootedness, where we refer to a polynomial as real-rooted if all of its roots are real. 
It is a well-known consequence of the general theory of real-rooted polynomials~\cite{StanleyLogConcave} that if $h^*(Z,x)$ has only real roots, its coefficient sequence is unimodal, but the converse does not hold.
Recent work by Beck, Jochemko, and McCullough has shown that $h^*$-polynomials for zonotopes are as well-behaved as possible from this perspective.

\begin{theorem}[Beck, Jochemko, and McCullough~\cite{beckjochemkomccullough}, Theorem 1.2]\label{thm:realrooted}
Let $Z$ be a $n$-dimensional lattice zonotope. 
Then the $h^*$-polynomial $h^*_Z(t)=h_0+h_1t+\dots+h_nt^n$ has only real roots. 
Moreover,
$$h_0\leq \cdots \leq h_{\frac{n}{2}}\geq \cdots \geq h_n\text{ if $n$ is even }$$
and 
$$h_0\leq \cdots \leq h_{\frac{n-1}{2}} \text{ and } h_{\frac{n+1}{2}}\geq \cdots \geq h_n \text{ if $n$ is odd.}$$
\end{theorem}

\begin{example}\label{ex:h*unimodal} 
Consider the zonotope $Z$ from Example~\ref{ex:stanleytheorem}.
We saw that $L_Z(t)=7t^2+4t+1$, and it is a nice exercise to change basis and show that $h^*(Z;x)=1+9x+4x^2$.
It follows easily that $h^*(Z;x)$ has only real roots and is clearly unimodal.
\end{example}

The key idea in the proof of Theorem~\ref{thm:realrooted} is again the half-open decomposition of $Z$.
The paper by Beck, Jochemko, and McCullough~\cite{beckjochemkomccullough} is well-written, and we recommend the interested reader look there for details.
However, we do want to discuss a connection with permutations and descents that is a core ingredient of their proof.

Let $\mathfrak{S}_n$ denote the set of all permutations on $[n]:=\{1,2,3,\ldots,n\}$. 
For a permutation word $\sigma=\sigma_1\sigma_2\cdots\sigma_n$ in $\mathfrak{S}_n$ the \emph{descent set} is defined by 
\begin{equation*}
\mathrm{Des}(\sigma):=\{i\in [n-1]: \sigma_i>\sigma_{i+1}\}\, . 
\end{equation*}
The number of descents of $\sigma$ is denoted by $\mathrm{des}(\sigma):=|\mathrm{Des}(\sigma)|$.
The \emph{Eulerian number} $a(n,k)$ counts the number of permutations in $\mathfrak{S}_n$ with exactly $k$ descents: 
\begin{equation*}
a(n,k):=|\{\sigma \in \mathfrak{S}_n: \mathrm{des}(\sigma)=k\}|\, . 
\end{equation*}
The $(A,j)$-\emph{Eulerian number}
\begin{equation*}
a_j(n,k):=|\{\sigma \in \mathfrak{S}_n: \sigma_n=n+1-j \text{ and } \mathrm{des}(\sigma)=k\}|
\end{equation*}
is a refinement of the descent statistic which gives the number of permutations $\sigma\in \mathfrak{S}_n$ with last letter $n+1-j$ and exactly $k$ descents. 
The associated polynomial is known as the $(A,j)$-\emph{Eulerian polynomial}, sometimes called the \emph{restricted Eulerian polynomial}, and is defined as 
\begin{equation*}
A_j(n,x):=\sum_{k=0}^{n-1}a_j(n,k)x^k \, .
\end{equation*}
The half-open decomposition of a lattice zonotope leads to the following result.

\begin{theorem}[Beck, Jochemko, and McCullough~\cite{beckjochemkomccullough}, Theorem 1.3]\label{thm:hstarcone}
Let $n\geq 1$. 
The convex hull of the $h^*$-polynomials of all $n$-dimensional lattice zonotopes (viewed as points in $\R^{n+1}$) and the convex hull of the $h^*$-polynomials of all $n$-dimensional lattice parallelepipeds are both equal to the $n$-dimensional simplicial cone 
$$A_1(n+1,x)+\R_{\geq 0}A_2(n+1,x)+\cdots+\R_{\geq 0}A_{n+1}(n+1,x).$$
\end{theorem}

It is fascinating that the convex hull of zonotope $h^*$-vectors admits such a beautiful combinatorial description, while the convex hull of the set of all $h^*$-vectors for lattice polytopes of a fixed dimension appears to be much more complicated~\cite{BDDPS,stapledoninequalities}.
Given our knowledge of Theorem~\ref{thm:hstarcone}, it would be interesting if the following problem turns out to be tractable, at least in low dimensions.

\begin{problem}
For fixed $n$, characterize the set of all $h^*$-vectors of $n$-dimensional parallelepipeds/zonotopes.
\end{problem}

\section{Graphical and Laplacian Zonotopes}\label{sec:graphical}

There are several interesting zonotopes associated with a finite simple graph $G$. 
In this section we discuss two such constructions, each of which is related to the number of spanning trees of $G$.
For the first construction, we recall that for a polynomial $f=\sum_{\ba\in \Z^n}\beta_{\ba} t_1^{a_1}\cdots t_n^{a_n}$, the \emph{Newton polytope} $\mathrm{Newton}(f)$  is the convex hull of integer points $\ba\in \Z^n$ such that $\beta_{\ba}\neq 0$. 
It is known that $\mathrm{Newton}(f\cdot g)$ is the Minkowski sum $\mathrm{Newton}(f)+\mathrm{Newton}(g)$, which is the fundamental ingredient for the following definition. 

\begin{definition}\label{def:graphicalz}
For a graph $G$ on the vertex set $[n]$, the \emph{graphical zonotope} $Z_G$ is defined to be
\begin{equation*}
Z_G:=\sum_{(i,j)\in G}[e_i,e_j]=\mathrm{Newton}\left(\prod_{(i,j)\in G}(t_i-t_j)\right) \, ,
\end{equation*}
where the Minkowski sum and the product are over edges $(i,j)$, $i<j$, of the graph $G$, and $e_1,\dots, e_n$ are the coordinate vectors in $\R^n$. 
\end{definition}

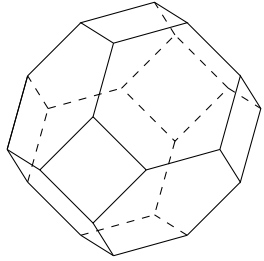
\begin{figure}
\begin{center}
\begin{tikzpicture}
\draw[black] (1.414,.707,0) -- (1.414,0,.707) -- (1.414,-.707,0) -- (1.414,0,-.707) -- cycle;
\draw[black] (-1.414,.707,0) -- (-1.414,0,.707) ;
\draw[black] (-1.414,.707,0) -- (-1.414,0,.707);
\draw[black] (-1.414,0,.707) -- (-1.414,-.707,0);
\draw[black,dashed]  (-1.414,-.707,0) -- (-1.414,0,-.707); 
\draw[black,dashed] (-1.414,0,-.707) -- (-1.414,.707,0);
\draw[black] (.707,1.414,0) -- (0,1.414,.707) -- (-.707,1.414,0) -- (0,1.414,-.707) -- cycle;
\draw[black]  (.707,-1.414,0) -- (0,-1.414,.707);
\draw[black] (0,-1.414,.707) -- (-.707,-1.414,0);
\draw[black,dashed] (-.707,-1.414,0) -- (0,-1.414,-.707);
\draw[black,dashed] (.707,-1.414,0) -- (0,-1.414,-.707) ;
\draw[black] (.707,0,1.414) -- (0,.707,1.414) -- (-.707,0,1.414) -- (0,-.707,1.414) -- cycle;
\draw[black,dashed] (.707,0,-1.414) -- (0,.707,-1.414) -- (-.707,0,-1.414) -- (0,-.707,-1.414) -- cycle;
\draw (1.414,.707,0) -- (.707,1.414,0)
  (1.414,0,.707) -- (.707,0,1.414)
  (0,1.414,.707) -- (0,.707,1.414)
  (1.414,-.707,0) -- (.707,-1.414,0)
  (-.707,0,1.414) -- (-1.414,0,.707)
  (0,-.707,1.414) -- (0,-1.414,.707)
  (-.707,1.414,0) -- (-1.414,.707,0)
  (-1.414,-.707,0) -- (-.707,-1.414,0);
  \draw[black,dashed] (-1.414,0,-.707) -- (-.707,0,-1.414);
  \draw[black,dashed] (0,.707,-1.414)  -- (0,1.414,-.707);
  \draw[black,dashed] (0,-.707,-1.414) -- (0,-1.414,-.707);
  \draw[black,dashed] (.707,0,-1.414) -- (1.414,0,-.707);
\end{tikzpicture}
\end{center}
\caption{The permutahedron $\mathcal{P}_4$.}
\label{fig:p4}
\end{figure}

\begin{example}\label{ex:perm}
The \emph{$n$-permutahedron} is the polytope in $\R^n$ whose vertices are the $n!$ permutations of $[n]$:
\begin{equation*}
\mathcal{P}_n:=\mathrm{conv}\{(\pi(1),\pi(2),\dots,\pi(n)): \pi\in \mathfrak{S}_n\}.
\end{equation*}
Using the Vandermonde determinant, it is straightforward to show that $\mathcal{P}_n$ is the graphical zonotope $Z_{K_n}$ for the complete graph $K_n$. 
Figure~\ref{fig:p4} shows $\mathcal{P}_4$.
\end{example}

The following beautiful theorem shows that the volume and lattice points of $Z_G$ encode the number of spanning trees and forests in $G$.

\begin{theorem}[Stanley~\cite{StanleyVol1}, Exercise 4.32; Postnikov~\cite{PostnikovAssPermBeyond}, Proposition 2.4]\label{thm:graphzonotopespan}
For a connected graph $G$ on $n$ vertices, the volume of the graphical zonotope $Z_G$ equals the number of spanning trees of $G$. 
The number of lattice points of $Z_G$ equals the number of forests in the graph $G$. 
\end{theorem}

\begin{example}
The number of spanning trees of the connected graph $K_n$ is $n^{n-2}$ and by Theorem \ref{thm:graphzonotopespan} it is also the volume of $\mathcal{P}_n$.
 Furthermore, the number of lattice points of $\mathcal{P}_n$ equals the number of forests on $n$ labeled vertices. 
\end{example}

\begin{remark}
Gruji\'c~\cite{GrujicFaceGraphZonotope} has shown that the $f$-polynomial of $Z_G$, which encodes the number of faces of $Z_G$ in each dimension, can be obtained as the principal specialization of the $q$-analogue of the chromatic symmetric function of $G$.
\end{remark}

Another zonotope associated with a finite simple graph arises from the Laplacian matrix, defined as follows.

\begin{definition}\label{def:graphlaplacian}
Let $G$ be an undirected graph with vertex set $V(G)=\{v_1,\dots, v_p\}$. 
The \emph{Laplacian matrix} $\LP:=\LP(G)$ is the $p\times p$ matrix
\[
\LP_{ij}=\begin{cases}
-1, &\text{ if $i\neq j$ and there is an edge between}\\
& \text{     vertices $v_i$ and $v_j$}\\
\mathrm{deg}(v_i), &\text { if $i=j$ },
\end{cases}
\]
where $\mathrm{deg}(v_i)$ denotes the degree (number of incident edges) of $v_i$. 
\end{definition}

The Laplacian matrix of a graph admits a standard factorization that is relevant to our discussion.
Let $G$ have an arbitrary edge orientation.
Let $\partial$ denote the $|V|\times |E|$ matrix where
\[\partial_{(v,e)} =
\begin{cases}
-1, &\text{ if }v\text{ is the negative endpoint of }e\\
1 & \text{ if }v\text{ is the positive endpoint of }e\\
0 &\text { otherwise }
\end{cases} \, .
\]
We call $\partial$ the \emph{signed incidence matrix} for $G$.
It is a straightforward exercise to verify that $\LP(G)=\partial\partial^T$.

A recent theorem due to Dall and Pfeifle connects the geometry of zonotopes arising from $\partial$ and $\LP(G)$.

\begin{theorem}[Dall and Pfeifle~\cite{dallpfeifle}, Theorem 17]\label{thm:laplacianvolume}
Let $\partial$ and $\LP(G)$ be the signed incidence matrix and Laplacian, respectively, for $G$.
Let $\overline{\LP}$ be the matrix obtained by taking any basis for $\Z^{|V|}\cap \mathrm{span}_{\R}(\LP(G))$ from among the columns of $\LP(G)$.
Then the volume of the zonotope defined by the columns of $\partial$ is equal to the volume of the zonotope defined by the columns of $\overline{\LP}$.
\end{theorem}

The reason why Theorem~\ref{thm:laplacianvolume} is interesting is that it allowed Dall and Pfeifle to give a polyhedral proof of the Matrix Tree Theorem, one of the most important theorems related to graph Laplacians, stated as follows.

\begin{theorem}\label{thm:matrixtree}
Let $G$ be a finite connected $p$-vertex graph without loops.
Let $1\leq i\leq p$, and let $\LP_0$ denote $\LP(G)$ with the $i$-th row and column removed. Let $\kappa(G)$ denote the number of spanning trees of $G$. 
Then 
\begin{equation*}
\kappa(G)=|\mathrm{det}(\LP_0)|.
\end{equation*}
\end{theorem}

\begin{remark}
Another interesting zonotope related to graphs is the zonotope of degree sequences of length $n$.
This family of zonotopes was studied by Stanley~\cite{stanleyzonotope} from an Ehrhart-theoretic perspective.
\end{remark}

\section{Fixed Subpolytopes of Permutahedra}\label{sec:fixed}

Because of the close connection between permutahedra (defined in Section~\ref{sec:graphical}) and symmetric groups (defined in Section~\ref{sec:ehrharthstar}), it is of interest to investigate the geometric properties of the action of the symmetric group on the permutahedron.
Specifically, the symmetric group $\mathfrak{S}_n$ acts on $\mathcal{P}_n\subset \R^n$ by permuting coordinates, i.e., $\sigma\in \mathfrak{S}_n$ acts on a point $\bx=(x_1,x_2,\dots,x_n)\in \mathcal{P}_n$, by  $\sigma\cdot \bx:=(x_{\sigma^{-1}(1)},x_{\sigma^{-1}(2)},\dots, x_{\sigma^{-1}(n)})$.
Recent research~\cite{ardilaschindlervindasmelendez} has focused on the fixed points of permutations under this action.

\begin{definition}\label{def:fixedsubpoly}
Let $\sigma\in\mathfrak{S}_n$.
The \emph{fixed subpolytope of the permutahedron} relative to $\sigma$ is
\begin{equation*}
\mathcal{P}^{\sigma}_n=\{\bx\in \mathcal{P}_n: \sigma\cdot \bx=\bx\} \, ,
\end{equation*}
the subpolytope of the permutahedron $\mathcal{P}_n$ fixed $\sigma$.
\end{definition}

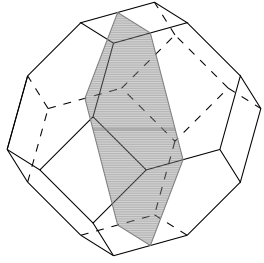
\begin{figure}
\begin{center}
\begin{tikzpicture}
\draw[black] (1.414,.707,0) -- (1.414,0,.707) -- (1.414,-.707,0) -- (1.414,0,-.707) -- cycle;
\draw[black] (-1.414,.707,0) -- (-1.414,0,.707) ;
\draw[black] (-1.414,.707,0) -- (-1.414,0,.707);
\draw[black] (-1.414,0,.707) -- (-1.414,-.707,0);
\draw[black,dashed]  (-1.414,-.707,0) -- (-1.414,0,-.707); 
\draw[black,dashed] (-1.414,0,-.707) -- (-1.414,.707,0);
\draw[black] (.707,1.414,0) -- (0,1.414,.707) -- (-.707,1.414,0) -- (0,1.414,-.707) -- cycle;
\draw[black]  (.707,-1.414,0) -- (0,-1.414,.707);
\draw[black] (0,-1.414,.707) -- (-.707,-1.414,0);
\draw[black,dashed] (-.707,-1.414,0) -- (0,-1.414,-.707);
\draw[black,dashed] (.707,-1.414,0) -- (0,-1.414,-.707) ;
\draw[black] (.707,0,1.414) -- (0,.707,1.414) -- (-.707,0,1.414) -- (0,-.707,1.414) -- cycle;
\draw[black,dashed] (.707,0,-1.414) -- (0,.707,-1.414) -- (-.707,0,-1.414) -- (0,-.707,-1.414) -- cycle;
\draw (1.414,.707,0) -- (.707,1.414,0)
  (1.414,0,.707) -- (.707,0,1.414)
  (0,1.414,.707) -- (0,.707,1.414)
  (1.414,-.707,0) -- (.707,-1.414,0)
  (-.707,0,1.414) -- (-1.414,0,.707)
  (0,-.707,1.414) -- (0,-1.414,.707)
  (-.707,1.414,0) -- (-1.414,.707,0)
  (-1.414,-.707,0) -- (-.707,-1.414,0);
  \draw[black,dashed] (-1.414,0,-.707) -- (-.707,0,-1.414);
  \draw[black,dashed] (0,.707,-1.414)  -- (0,1.414,-.707);
  \draw[black,dashed] (0,-.707,-1.414) -- (0,-1.414,-.707);
  \draw[black,dashed] (.707,0,-1.414) -- (1.414,0,-.707);

\coordinate (A) at (.3535,1.414,.3535);
\coordinate (B) at (-.3535,1.414,-.3535);
\coordinate (C) at (1.0605,0,1.0605);
\coordinate (D) at (.3535,-1.414,.3535);
\coordinate (E) at (-.3535,-1.414,-.3535);
\coordinate (F) at (-1.0605,0,-1.0605);
\draw[gray,top color=gray, bottom color=gray, fill opacity=0.35] (F) -- (E) -- (D) -- (C) -- (A) -- (B) -- cycle;
\end{tikzpicture}
\end{center}
\caption{The subpolytope $\mathcal{P}_4^{(12)}$ of the permutahedron $\mathcal{P}_4$ fixed by $(12)\in \mathfrak{S}_4$.}
\label{fig:fixedsubpolytope}
\end{figure}

Given the combinatorial nature of the normalized volume of $\mathcal{P}_n$ described in Example~\ref{ex:perm}, it is of interest to compute the normalized volumes of fixed subpolytopes of $\mathcal{P}_n$.
Before doing so, we first describe the vertices, supporting hyperplanes, and Minkowski sum representation of $\mathcal{P}^{\sigma}_n$.
Let $e_i$ denote the $i$-th standard basis vector in $\R^n$.

\begin{theorem}[Ardila, Schindler, Vindas-Mel\'endez~\cite{ardilaschindlervindasmelendez}, Theorem 2.12] \label{thm:fixedsubpolytope}
Let $\sigma$ be a permutation of $[n]$ having disjoint cycles $\sigma_1, \ldots, \sigma_m$, which have respective lengths $l_1, \ldots, l_m$. 
The fixed subpolytope $\mathcal{P}_n^\sigma$ can be described in the following four ways:
\begin{enumerate}
\item[0.]
It is the set of points $x$ in the permutahedron $\mathcal{P}_n$ such that $\sigma \cdot \bx = \bx$.
\item[1.]
It is the set of points $x \in \R^n$ satisfying
\begin{enumerate}
\item $x_1 + x_2 + \cdots + x_n = 1 + 2 + \cdots +  n$,
\item for any proper subset $\{i_1, i_2, \dots, i_k\} \subset \{1,2, \dots, n\}$, \[x_{i_1} + x_{i_2} + \cdots + x_{i_k} \leq 1 + 2 + \cdots + k, \textrm{ and} \]
\item
for any $i$ and $j$ which are in the same cycle of $\sigma$, $x_i=x_j$.
\end{enumerate}
\item[2.]
It is the convex hull of the set of vertices $\mathrm{Vert}(\sigma)$  consisting of the $m!$ points 
\begin{equation*}
\overline{\bov_{\prec}} := \sum_{k=1}^m \bigg( \dfrac{l_k +1}{2}
+ \sum_{j \, : \, \sigma_j \prec \sigma_k} l_j  \bigg) e_{\sigma_k}
\end{equation*}
as $\prec$ ranges over the $m!$ possible linear orderings of the disjoint cycles $\sigma_1, \sigma_2, \dots, \sigma_m$.
\item[3.]
It is the Minkowski sum 
\begin{equation*}
M_{\sigma}:=\sum_{1 \leq j < k \leq m}[l_je_{\sigma_k}, l_ke_{\sigma_j}]+ \sum_{k=1}^m \dfrac{l_k+1}{2}e_{\sigma_k}.
\end{equation*}
\end{enumerate}
Consequently, the fixed polytope $\mathcal{P}_n^\sigma$ is an $(m-1)$-dimensional zonotope that is combinatorially isomorphic to the permutahedron $\mathcal{P}_m$.
\end{theorem}

\begin{example}
In Figure \ref{fig:fixedsubpolytope}, we see that $\mathcal{P}_4^{(12)}$ has six vertices and Theorem \ref{thm:fixedsubpolytope} tells us that the vertices correspond to the 6 orderings of $\sigma_1=(12), \sigma_2=(3), \sigma_3=(4)$:
\begin{align*}
&\tfrac{1+2}{2}\, e_{12} +3e_{3} +4e_{4}, \quad
\tfrac{1+2}{2}\,e_{12} +4e_{3} +3e_{4},\\
&\tfrac{2+3}{2}\,e_{12} +1e_{3} +4e_{4}, \quad
\tfrac{2+3}{2}\,e_{12} +4e_{3} +1e_{4}, \\
&\tfrac{3+4}{2}\,e_{12} +1e_{3} +2e_{4}, \quad
\tfrac{3+4}{2}\,e_{12} +2e_{3} +1e_{4}.
\end{align*}
Furthermore, $\mathcal{P}_4^{(12)}$ is the Minkowski sum $$[e_{12},2e_{3}]+[e_{12},2e_{4}]+[e_3,e_4]+\frac{3}{2}e_{12}+e_3+e_4.$$
\end{example}

Using Theorem~\ref{thm:fixedsubpolytope}, it is possible to determine an arithmetical formula for the normalized volume of $\mathcal{P}_n$ related to only the disjoint cycle decomposition of $\sigma$.

\begin{theorem}[Ardila, Schindler, Vindas-Mel\'endez~\cite{ardilaschindlervindasmelendez}, Theorem 1.2]\label{thm:volfixed}
If $\sigma$ is a permutation of $[n]$ whose cycles have lengths $l_1, \ldots, l_m$, then the normalized volume of the fixed subpolytope of $\mathcal{P}_n$ is 
\begin{equation*}
\Vol\,  \mathcal{P}_n^\sigma = n^{m-2} \gcd(l_1, \ldots, l_m).
\end{equation*}
\end{theorem}

\section{Cyclic Zonotopes}\label{sec:cyclic}

Let $\{(1,t,t^2,t^3,\ldots,t^{n-1}):t\in \R\}$ be the \emph{moment curve} in $\R^n$.
The moment curve is often used in discrete geometry to select points in general position, and the well-studied family of cyclic polytopes arise by taking the convex hull of a finite set of points on the moment curve.
There is a zonotopal analogue of cyclic polytopes, defined as follows using centrally symmetric zonotopes.

\begin{definition}
The \emph{cyclic zonotope} $Z(d,n)$ defined by $t_1<t_2<\cdots<t_d$ is the $n$-zonotope $Z_0(\bov_1,\ldots,\bov_d)=Z(\pm \bov_1,\ldots,\pm \bov_d)$ generated by the $d$ vectors $\bov_i=(1,t_i,\dots,t_i^{n-1})$ in $\R^n$.
\end{definition}

It is known that while $Z(d,n)$ depends on the set of parameters $\{t_1,\dots,t_d\}$, its face poset depends only on $d$ and $n$, hence the suppression of the values of $t$ in our notation. 
When $\{t_1,\dots,t_d\}\subset \Z$, we say that $Z(d,n)$ is an \emph{integral} cyclic zonotope, though our discussion in this section will apply to cyclic zonotopes more generally.

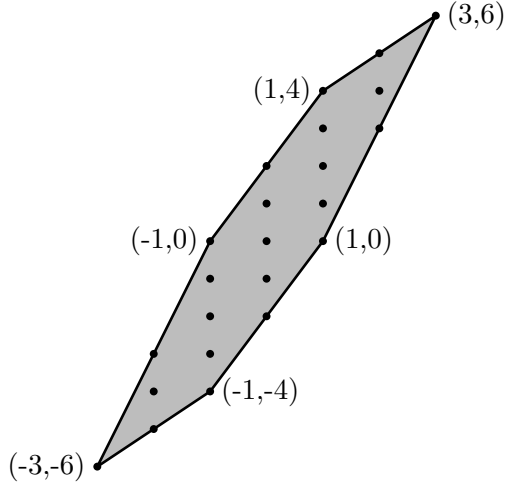
\begin{figure}
\begin{center}
\begin{tikzpicture}[x  = {(1.5cm,0cm)},
                    y  = {(0cm,1cm)},
                    z  = {(0cm,0cm)},
                    scale = .5,
                    color = {lightgray}]

  \definecolor{pointcolor_my_polytope}{rgb}{ 0,0,0 }
  \tikzstyle{pointstyle_my_polytope} = [fill=pointcolor_my_polytope]

  \coordinate (v0_my_polytope) at (-3, -6);
  \coordinate (v1_my_polytope) at (-1, -4);
  \coordinate (v2_my_polytope) at (1, 0);
  \coordinate (v3_my_polytope) at (3, 6);
  \coordinate (v4_my_polytope) at (1, 4);
  \coordinate (v5_my_polytope) at (-1, 0);

  \definecolor{edgecolor_my_polytope}{rgb}{ 0,0,0 }

  \definecolor{facetcolor_my_polytope}{rgb}{0.7, 0.7, 0.7}

  \tikzstyle{facestyle_my_polytope} = [fill=facetcolor_my_polytope, fill opacity=0.85, draw=edgecolor_my_polytope, line width=1 pt, line cap=round, line join=round]

  \draw[facestyle_my_polytope] (v1_my_polytope) -- (v2_my_polytope) -- (v3_my_polytope) -- (v4_my_polytope) -- (v5_my_polytope) -- (v0_my_polytope) -- (v1_my_polytope) -- cycle;

  \node at (-3,-6) [text=black, inner sep=4pt, left, draw=none, align=left] {(-3,-6)};
  \node at (3,6) [text=black, inner sep=4pt, right, draw=none, align=left] {(3,6)};
  \node at (-1,-4) [text=black, inner sep=4pt, right, draw=none, align=left] {(-1,-4)};
  \node at (-1,0) [text=black, inner sep=4pt, left, draw=none, align=left] {(-1,0)};
  \node at (1,0) [text=black, inner sep=4pt, right, draw=none, align=left] {(1,0)};
  \node at (1,4) [text=black, inner sep=4pt, left, draw=none, align=left] {(1,4)};

  \fill[pointcolor_my_polytope] (-3,-6) circle (3 pt);
  \fill[pointcolor_my_polytope] (-2,-3) circle (3 pt);
  \fill[pointcolor_my_polytope] (-2,-4) circle (3 pt);
  \fill[pointcolor_my_polytope] (-2,-5) circle (3 pt);
  \fill[pointcolor_my_polytope] (-1,-4) circle (3 pt);
  \fill[pointcolor_my_polytope] (-1,-3) circle (3 pt);
  \fill[pointcolor_my_polytope] (-1,-2) circle (3 pt);
  \fill[pointcolor_my_polytope] (-1,-1) circle (3 pt);
  \fill[pointcolor_my_polytope] (-1,0) circle (3 pt);
  \fill[pointcolor_my_polytope] (0,2) circle (3 pt);
  \fill[pointcolor_my_polytope] (0,1) circle (3 pt);
  \fill[pointcolor_my_polytope] (0,0) circle (3 pt);
  \fill[pointcolor_my_polytope] (0,-1) circle (3 pt);
  \fill[pointcolor_my_polytope] (0,-2) circle (3 pt);
  \fill[pointcolor_my_polytope] (1,0) circle (3 pt);
  \fill[pointcolor_my_polytope] (1,1) circle (3 pt);
  \fill[pointcolor_my_polytope] (1,2) circle (3 pt);
  \fill[pointcolor_my_polytope] (1,3) circle (3 pt);
  \fill[pointcolor_my_polytope] (1,4) circle (3 pt);
  \fill[pointcolor_my_polytope] (2,3) circle (3 pt);
  \fill[pointcolor_my_polytope] (2,4) circle (3 pt);
  \fill[pointcolor_my_polytope] (2,5) circle (3 pt);
  \fill[pointcolor_my_polytope] (3,6) circle (3 pt);

\end{tikzpicture}
\end{center}
\caption{The cyclic zonotope $Z(3,2)$ generated by $t\in\{1,2,3\}$.}
\label{fig:cycliczonotope}
\end{figure}

\begin{example}
Figure \ref{fig:cycliczonotope} illustrates the cyclic zonotope $Z(3,2)$ generated by $t\in\{1,2,3\}$.
\end{example}

Cyclic zonotopes have particularly nice geometric properties, especially with regard to subdivisions.
We have seen throughout this paper the importance of the half-open decompositions of zonotopes established in Theorem~\ref{thm:decomp}.
A relaxation of this idea is the notion of a zonotopal subdivision, which we now review in the setting of $Z(d,n)$.
Denote by $\Lambda_d$ the set $\{-,0,+\}^d$ of sign vectors of length $d$ and write $X=(X_1,X_2,\dots, X_d)$ for $X\in \Lambda_d$. 
The set $\Lambda_d$ is partially ordered by extending coordinatewise the partial order on $\{-,0,+\}$ defined by the relations $0<-$ and $0<+$. 
A \emph{subzonotope} $Z_X$ of $Z(d,n)$ is a pair $(X,Z_X)$ of a sign vector $X\in \Lambda_d$ and the translated zonotope
\begin{equation*}
Z_X=\sum_{X_i=+}\bov_i- \sum_{X_i=-}\bov_i+\sum_{X_i=0}[-\bov_i,\bov_i] \, .
\end{equation*}
We say that $Z_C$ is a \emph{face} of $Z_X$ if there exists a face $F$ of the polytope $Z_X$ such that $C$ is minimum among all sign vectors $Y\geq X$ with $Z_Y=F$. 
Two subzonotopes $Z_X$ and $Z_Y$ of $Z(d,n)$ are said to \emph{intersect properly} if either $Z_X \cap Z_Y = \emptyset$ or there exists a sign vector $C$ such that 
\begin{enumerate}
\item $Z_X\cap Z_Y=Z_C$ and
\item $Z_C$ is a face of both $Z_X$ and $Z_Y$.
\end{enumerate}


\begin{definition}
A \emph{zonotopal subdivision} $S$ of $Z(d,n)$ is a set of subzonotopes of $Z(d,n)$ with the following properties: 
\begin{enumerate}
\item if $Z_C$ is a face of $Z_X\in S$ then $Z_C\in S$,
\item any two elements of $S$ intersect properly and 
\item $\bigcup_{Z_X\in S} Z_X=Z(d,n)$.
\end{enumerate}
\end{definition}


As with triangulations of polytopes, there is a notion of shellability for zonotopal subdivisions.

\begin{definition}
A zonotopal subdivision $S$ of $Z(d,n)$ is \emph{shellable} if it is pure of dimension $n$, i.e., each subzonotope of $Z(d,n)$ has dimension $n$, and there exists a linear ordering of the subzonotopes $Z_{{X_1}},\dots,Z_{{X_k}}$, called a \emph{shelling}, such that for $1<j\leq k$, 
\[
Z_{X_j}\cap \left(\bigcup_{i\leq i<j}Z_{X_i}\right)
\]
is not empty and equals the union of the $(n-1)$-polytopes in an initial segment of a shelling of the boundary of $Z_{X_j}$.
\end{definition}

\begin{figure}
\begin{center}
\begin{tikzpicture}[scale=1]
\draw[gray,  thick, gray,top color=gray, bottom color=white, fill opacity=0.4] (7,1) -- (8,1) -- (8,0) -- (7,0)--cycle;
\draw[gray,  thick, gray,top color=gray, bottom color=white, fill opacity=0.4] (8,1) -- (9,2)-- (8,2) --(7,1) -- cycle;
\draw[gray,  thick, gray,top color=white, bottom color=gray, fill opacity=0.4](8,0) -- (9,1)--(9,2)--(8,1)--cycle;
\draw[gray,  thick, gray,top color=white, bottom color=gray, fill opacity=0.3](8,0) -- (9,1)--(10,0)--(9,-1)--cycle;
\draw[gray,  thick, gray,top color=gray, bottom color=white, fill opacity=0.4] (8,0) -- (9,-1)--(8,-1)--(7,0)--cycle;
\draw[gray,  thick, gray,top color=gray, bottom color=white, fill opacity=0.2] (9,1) -- (9,2) -- (10,1) -- (10,0) --cycle;
\draw[black, ultra thick] (7,0) -- (7,1) -- (8,2) -- (9,2) -- (10,1) -- (10,0) -- (9,-1) -- (8, -1) -- cycle; 

 \node at (7,0) [text=black, inner sep=10pt, above right, draw=none, align=left] {1};
 \node at (8,-1) [text=black, inner sep=10pt, above, draw=none, align=left] {2};
 \node at (8,1) [text=black, inner sep=10pt, above, draw=none, align=left] {3};
 \node at (8,1) [text=black, inner sep=10pt, right, draw=none, align=left] {4};
 \node at (9,-1) [text=black, inner sep=25pt, above, draw=none, align=left] {5};
 \node at (9,1) [text=black, inner sep=10pt, right, draw=none, align=left] {6};

\end{tikzpicture}
\end{center}
\caption{A zonotopal subdivision of $Z_0((0,1),(1,0),(1,1),(1,-1))$ with a shelling.}
\label{fig:zonotopalsubdivision}
\end{figure}
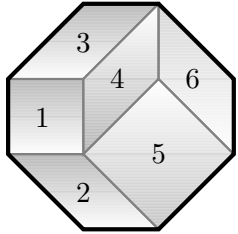

\begin{example}
Consider the zonotope $Z_0((0,1),(1,0),(1,1),(1,-1))$ illustrated in Figure~\ref{fig:fourlinesegmentzonotope}. The zonotope is shellable and Figure \ref{fig:zonotopalsubdivision} shows a shelling, that is, the subdivision given by the numbered facets together with their faces is shellable. 
\end{example}

A lovely result about cyclic zonotopes is the following theorem due to Athanasiadis.

\begin{theorem}[Athanasiadis~\cite{athanasiadiszonotope}, Theorem 1.2]\label{thm:cyclicshellable}
All zonotopal subdivisions of a cyclic zonotope are shellable. 
\end{theorem}

In addition to Theorem~\ref{thm:cyclicshellable}, Athanasiadis~\cite{athanasiadiszonotope} also proved that posets of proper zonotopal subdivisions of $Z(d,n)$ induced by canonical projections of cyclic zonotopes are homotopy equivalent to spheres.
Thus, there are interesting connections between the structure of cyclic zonotopes and poset topology.


\section*{Acknowledgments}

This work was partially supported by NSF Graduate Research Fellowship DGE-1247392 (ARVM).
Thanks to Federico Ardila, Matt Beck, and Sam Hopkins for helpful comments.



\bibliographystyle{ws-procs961x669}
\bibliography{Braun}

\end{document}